\documentclass[12pt]{amsart}
\usepackage{graphicx}
\usepackage{amsmath}
\usepackage{amsfonts}
\usepackage{amssymb}
\usepackage{amscd}
\usepackage{amsthm}

\newtheorem{theorem}{Theorem}
\newtheorem{corollary}[theorem]{Corollary}
\newtheorem{lemma}[theorem]{Lemma}
\newtheorem{proposition}[theorem]{Proposition}

\newtheorem{remark}[theorem]{Remark}

\def\X{\mathcal{E}}

\def\b{\beta}

\def\S{\mathbf{S}}

\def\MCG{\textup{MCG}}
\def\L{\mathcal{L}}
\def\PMCG{\textup{PMCG}}
\def\H{\textup{H}}
\def\T{\textup{T}}
\def\V{V}
\def\E{\bar{\mathcal{E}}}

\def\P{\mathcal{P}}
\date{\today}

\begin{document}

\title[Extending homeomorphisms from punctured surfaces]{Extending homeomorphisms from punctured surfaces to handlebodies}

\author{Alessia Cattabriga \and  Michele Mulazzani}

\begin{abstract}
Let $\H_g$ be a genus $g$ handlebody and $\MCG_{2n}(\T_g)$  be the
group of the isotopy classes of orientation preserving
homeomorphisms of $\T_g=\partial\H_g$,  fixing a given set of $2n$
points. In this paper we find a finite set of generators for
$\X_{2n}^g$, the subgroup of $\MCG_{2n}(\T_g)$ consisting of the
isotopy classes of homeomorphisms of $\T_g$  admitting an
extension to the handlebody and keeping fixed the union of  $n$
disjoint properly embedded trivial arcs. This result generalizes a
previous one obtained by the authors for  $n=1$. The subgroup
$\X_{2n}^g$ turns out to be important for the study of knots and
links  in closed 3-manifolds via $(g,n)$-decompositions. In fact,
the links represented by the isotopy classes belonging to the same
left cosets of $\X_{2n}^g$ in $\MCG_{2n}(\T_g)$ are equivalent.

\bigskip

\noindent{{\it Mathematics Subject
Classification 2000:} Primary  20F38; Secondary  57M25.}

\end{abstract}

\maketitle

\section{Introduction and preliminaries} \label{intro}

Let $\H_g$ be an oriented handlebody of genus $g\geq 0$ and
$\partial \H_g=\T_g$. Consider a system of $n$ disjoint properly
embedded trivial arcs\footnote{A set of mutually disjoint arcs
$\{A_1,\ldots ,A_n\}$ properly embedded in a handlebody $\H_g$ is
trivial if there exist $n$ mutually disjoint embedded discs,
called \textit{trivializing discs}, $D_1,\ldots ,D_n\subset \H_g$
such that $A_i\cap D_i= A_i\cap\partial D_i=A_i$, $A_i\cap
D_j=\emptyset$ and $\partial D_i-A_i\subset\partial \H_g$ for all
$i,j=1,\ldots ,n$ and $i\neq j$.} $\mathcal A=\{A_1,\ldots,A_n\}$
in $\H_g$ and let $P_{i1},P_{i2}$ be the endpoints of the arc
$A_i$, for $i=1,\ldots,n$. We denote with $\MCG_{2n}(\T_g)$ (resp.
$\MCG_n(\H_g)$),  the group of the isotopy classes of orientation
preserving homeomorphisms of $\T_g$ (resp. $\H_g$)  fixing the set
$\P_n=\{P_{i1},P_{i2}\,\vert\, i=1,\ldots,n\}$ (resp.
$A_1\cup\cdots \cup A_n$). The group $\MCG_{2n}(\T_g)$ is widely
studied and different finite presentations of it are known (see
 \cite{G,LP}). In this article we are interested in
studying the subgroup $\X_{2n}^g$ of $\MCG_{2n}(\T_g)$, which is
the image of the homomorphism $\MCG_n(\H_g)\to \MCG_{2n}(\T_g)$
induced by restriction.  In other words, an element
$f\in\MCG_{2n}(\T_g)$ belongs to $\X_{2n}^{g}$ if it admits an
extension to $\H_g$ fixing $A_1\cup\cdots \cup A_n$. If $h\in
\X_{2n}^g$ we still denote with $h$ an extension of it. The case
$n=0$  is studied in \cite{Su},  while in \cite{CM} we study the
case  $n=1$. Moreover, in \cite{H}, the case $g=0$ is
investigated. In this paper we find a finite set of generators for
$\X_{2n}^g$ and describe their extension to $\H_g$, for each
$n\geq 1$ and $g\geq 0$.

The main motivation for studying such subgroups lies in their
importance in the representation of knots and links via
$(g,n)$-decompositions. The notion of $(g,n)$-decomposition for
links in orientable closed connected 3-manifolds, given by Doll in
\cite{Do}, extends the one of bridge (or plat) decomposition for
links in $\S^3$.  Roughly speaking, a $(g,n)$-decomposition of a
link $L\subset M$ is the data of a Heegaard surface of genus $g$
in $M$ which cuts the link into two set of $n$ disjoint trivial
arcs (the underpasses and the overpasses). It is easy to see that
each link admits a $(g,n)$-decomposition, for suitable $g$  and
$n$ with $g\geq g_M$ and $n\geq n_L$, where $g_M$ denotes the
Heegaard genus of $M$ and $n_L$ is the number of components of
$L$.

The use of $(g,1)$-decompositions of knots is revealed to be very
fruitful in order to study different topics.  For example  the
strongly-cyclic branched coverings of knots in 3-manifolds are
analyzed in \cite{CMV}, while the Alexander polynomial of knots in
rational homology spheres is investigated in \cite{K}. The case of
knots admitting $(1,1)$-decompositions has been widely studied
(see for example \cite{CM1,CM2,CK,E,GHS,Ha,Sa}). In  \cite{GMM},
the Heegaard Floer homology of $(1,1)$-knots is computed, while in
\cite{S}, the geometry of such   knots is determined in terms of
the distance of the curve complex associated with their
$(1,1)$-decompositions.

In Section \ref{decomposition} we give the definition of
$(g,n)$-decomposition of links and describe its connections with
$\X_{2n}^g$. In Section  \ref{zero} we discuss the genus zero
case, while
 the general case is analyzed in
Section \ref{general}.

Now we briefly recall the definition of spin of a point along a
curve (see \cite{Bi}). Let $P$ be a point on $\T_g$ and $c$ be a
simple closed oriented curve on $\T_g$, containing $P$. Consider a
neighborhood $N$ of $c$, parametrized by coordinates
$(y,\theta)\in [-1,1]\times \S^1$, where $c$ is defined by $y=0$
and  $P$ is the point of coordinates $(0,0)$. The \textit{spin} of
$P$ about $c$ is the homeomorphism $s_{P,c}$ of $\T_g$ obtained
extending with the identity  the map defined on $N$ by
$(y,\theta)\mapsto(y,\theta+\theta')$ where
$$\theta'=\left\{  \begin{array}{ll} 2\pi (2y+1)& \textup{ if } -1\le y\le
-1/2\\0 &\textup{ if } -1/2\le y\le 1/2\\-2\pi (2y-1)& \textup{ if
} 1/2\le y\le 1
\end{array}\right..$$ If we
denote with $t_{\delta}$ the right-handed Dehn twist along the
simple closed curve $\delta$, then $s_{P,c}=t_{c_1}^{-1}t_{c_2}$,
where $c_1$ and $c_2$ are the curves corresponding to $y=-3/4$ and
$y=3/4$ respectively.

\section{$(g,n)$-decompositions of links}
\label{decomposition}

In  this section we recall the notion of $(g,n)$-decompositions of
links and describe the connection with $\X_{2n}^g$.

\begin{figure}[ht]
\begin{center}
\includegraphics*[totalheight=5.5cm]{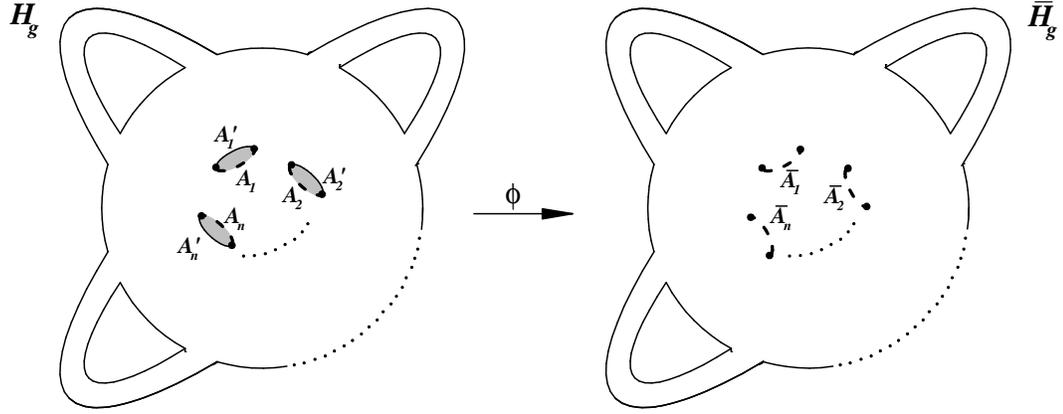}
\end{center}
\caption{A $(g,n)$-decomposition.} \label{deco}
\end{figure}

A $(g,n)$-decomposition for a link $L$ in a orientable closed
connected 3-manifold $M$ is the data
$$(M,L)=(\H_g, A_1\cup\cdots\cup A_n)\cup_{\phi}(\bar \H_g, \bar A_1\cup\cdots\cup \bar A_n)$$
where $\H_g,\bar \H_g$ are two oriented handlebodies of genus $g$,
$\mathcal A=\{A_1,\ldots,A_n\},\bar{\mathcal A}=\{\bar
A_1,\ldots,\bar A_n\}$ are two systems of $n$ properly embedded
trivial arcs in $\H_g$ and $\bar\H_g$ respectively, and
$\phi:(\partial\H_g,
\partial A_1\cup\cdots\cup \partial A_n)\to (\partial\bar\H_g,
\partial \bar A_1\cup\cdots\cup \partial \bar A_n)$ is an
 attaching (orientation preserving) homeomorphism (see Figure
\ref{deco}).

\begin{figure}[ht]
\begin{center}
\includegraphics*[totalheight=8cm]{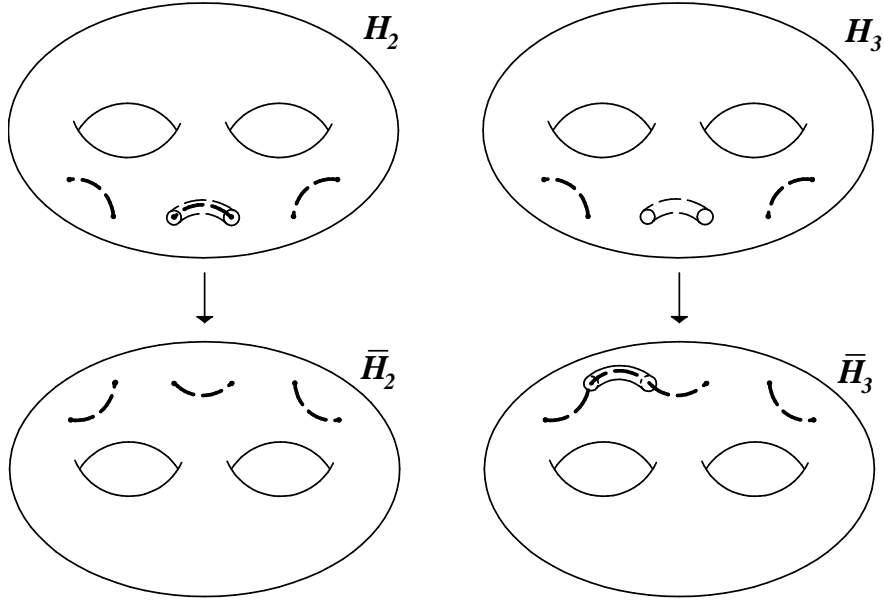}
\end{center}
\caption{From a $(2,3)$-decomposition to a $(3,2)$-decomposition.}
\label{components}
\end{figure}

A link admitting a $(g,n)$-decomposition will be called a
$(g,n)$-link. It is easy to see  that $L$ admits  a
$(g_M,n)$-decomposition, for a suitable $n$, where $g_M$ denotes
the Heegaard genus of $M$. Indeed,  given a Heegaard surface $S$
of genus $g_M$ for $M$,   there exists an immersion of $L$ in $S$
with a finite number of singular points which are double points.
So, by a slight modification of this immersion near each double
point, we are able  to embed $L$ in $S$, except for a finite
trivial set of arcs.  Moreover, if $n_L$ denotes the number of
components of $L$, by choosing a  sufficiently large $g$, it is
always possible to find a $(g,n_L)$-decomposition of $L$, since if
$n> n_L$ a $(g,n)$-decomposition determines a
$(g+1,n-1)$-decomposition (see Figure \ref{components}).

Let $L\subset M$ be a $(g,n)$-link with $(g,n)$-decomposition
$(M,L)=(\H_g, A_1\cup\cdots\cup A_n)\cup_{\phi}(\bar \H_g, \bar
A_1\cup\cdots\cup \bar A_n)$ and  let $\Upsilon:\bar\H_g \to \H_g$
be a fixed (orientation preserving) homeomorphism such that
$\Upsilon(\bar A_i)= A_i$, for $i=1,\ldots,n$, then
$\varphi=\Upsilon_{|\partial \bar H_g} \phi$ is an orientation
preserving homeomorphism of $(\partial \H_g,\partial
A_1\cup\cdots\cup\partial A_n)$.

Moreover, since two isotopic attaching homeomorphisms produce
equivalent $(g,n)$-links, we have a natural surjective map from
$\MCG_{2n}(\T_g)$, the  mapping class group of the $2n$-punctured
surface of genus $g$, to the class $\L_{g,n}$ of all $(g,n)$-links
$$\Theta_{g,n}:\MCG_{2n}(\T_g)\to\L_{g,n}$$
sending $\varphi\in \MCG_{2n}(\T_g)$ to the link $L_{\varphi}\in
\L_{g,n}$, associated to the attaching homeomorphism
$\Upsilon_{|\partial  H_g}^{-1}\varphi$. Unfortunately
$\Theta_{g,n}$ is far away from being  injective. Indeed, if we
denote with $\X_{2n}^g$  the subgroup of $\MCG_{2n}(\T_g)$
consisting of the elements that admit an extension to $\H_g$,
fixing $A_1\cup\cdots\cup A_n$ as a set, then for each
$\varepsilon\in\X_{2n}^g$, the link $L_{\varepsilon}$ is the
$n$-component trivial link  in the connected sum of $g$ copies of
$\S^2\times\S^1$. Moreover, for each $\varphi\in\MCG_{2n}(\T_g)$,
the links $L_{\varphi\varepsilon}$ and $L_{\varphi}$ are
equivalent. So we can restrict ourselves to considering the
surjective map
$$\Theta'_{g,n}:\MCG_{2n}(\T_g)/\X_{2n}^g\to \L_{g,n},$$
where $\MCG_{2n}(\T_g)/\X_{2n}^g$ denotes the set of left cosets
of $\X_{2n}^g$ in $\MCG_{2n}(\T_g)$. So, in this contest, it is
important to obtain information about $\X_{2n}^g$. In this paper
we find a finite set of generators for this  group.  In order to
do this, let us describe some elements of $\X_{2n}^g$.
\begin{description}

\item[Intervals] For $i=1,\ldots,n$, let  $D_i$ be  a trivializing
disc for $A_i$. Consider a tubular neighborhood $N$ of $D_i$ such
that $D_j\cap N=\emptyset$ for $j\ne i$ and parametrize it by
$N\cong \mathbf{D}^2\times [0,1]$ as in Figure \ref{finale}. For
$i=1,\ldots,n$, we denote with $\iota_i$ the homeomorphism of
$\H_g$ obtained extending by the identity the one defined on $N$
by $(\rho,\theta,t)\mapsto (\rho,\theta+\theta',t)$ where

\begin{equation}
\label{formula} \theta'=\left\{\begin{array}{ll}\pi &\textup{if }
0\le \rho,t\le 1/2
\\ 2\pi(1-\rho)
&\textup{if } 1/2\le \rho\le 1,\ 0\le t\le 1/2 \\
2\pi(1-t) &\textup{if } 0\le \rho \le 1/2,\ 1/2\le t\le 1
\\  4\pi(1-t)(1-\rho) & \textup{if } 1/2\le \rho,t \le 1.\end{array}\right.
\end{equation}

By definition, $\iota_i$ exchanges the endpoints of the arc $A_i$
and fixes pointwise the other arcs.

\begin{figure}[ht]
\begin{center}
\includegraphics*[totalheight=4cm]{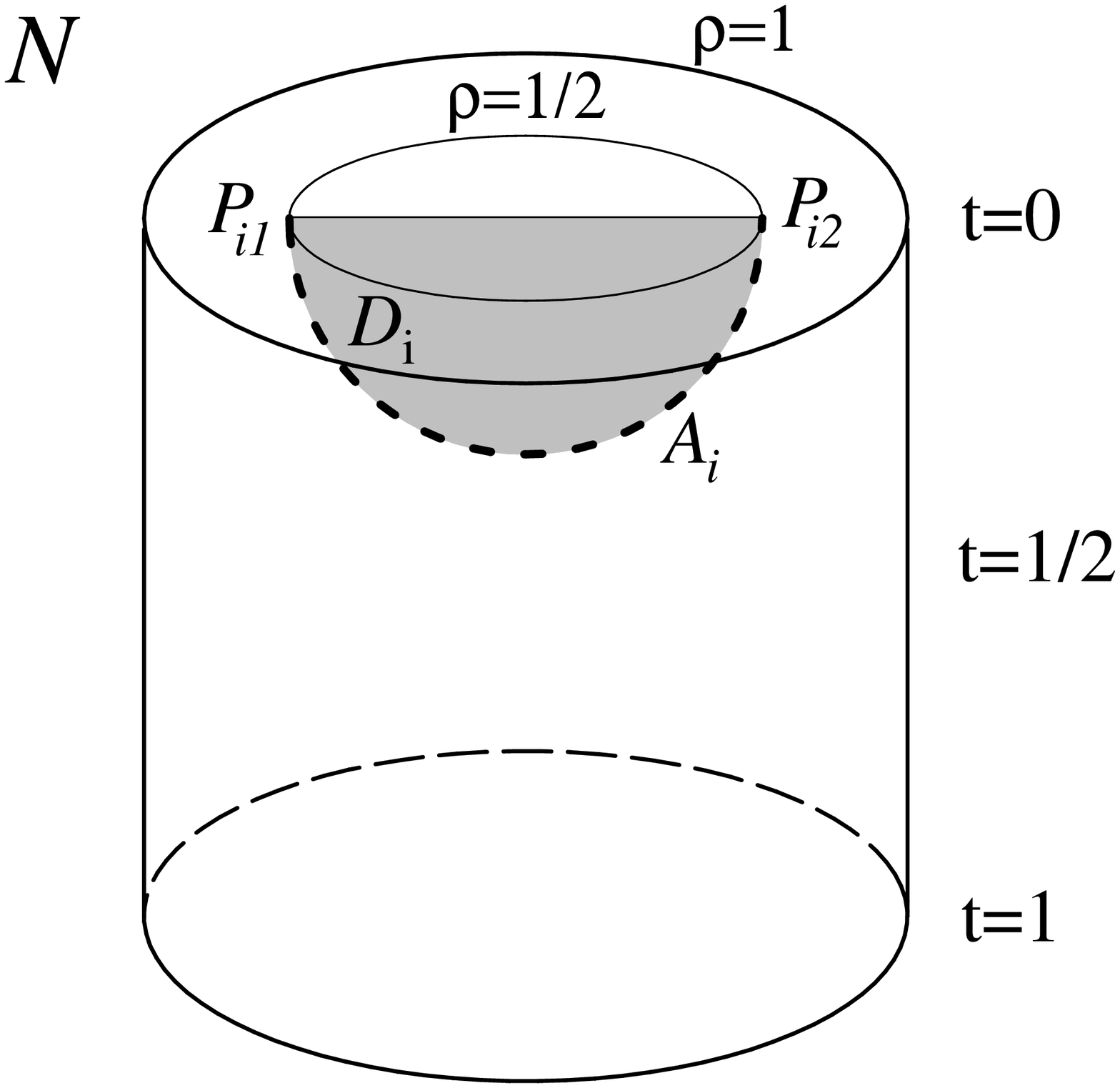}
\end{center}
\caption{} \label{finale}
\end{figure}

\item[Exchanging two arcs] Consider a tubular neighborhood $N$ of
$D_i\cup B \cup D_{i+1}$, where  $B$ is a band connecting $D_i$
and $D_{i+1}$, such that $D_j\cap N=\emptyset$, for $j\ne i,i+1$,
and parametrize it by $N\cong \mathbf{D}^2\times [0,1]$ as in
Figure \ref{finale2}. For $i=1,\ldots,n-1$, we denote with
$\bar{\lambda}_i$ the homeomorphism of $\H_g$ defined as in
(\ref{formula}) on $N$ and as the identity outside $N$. Set
$\lambda_i=\iota_{i+1}^{-1}\iota_i^{-1}\bar{\lambda}_i$ exchanges
the arcs $A_i$ and $A_{i+1}$ and fixes pointwise the other arcs.
\begin{figure}
\begin{center}
\includegraphics*[totalheight=5cm]{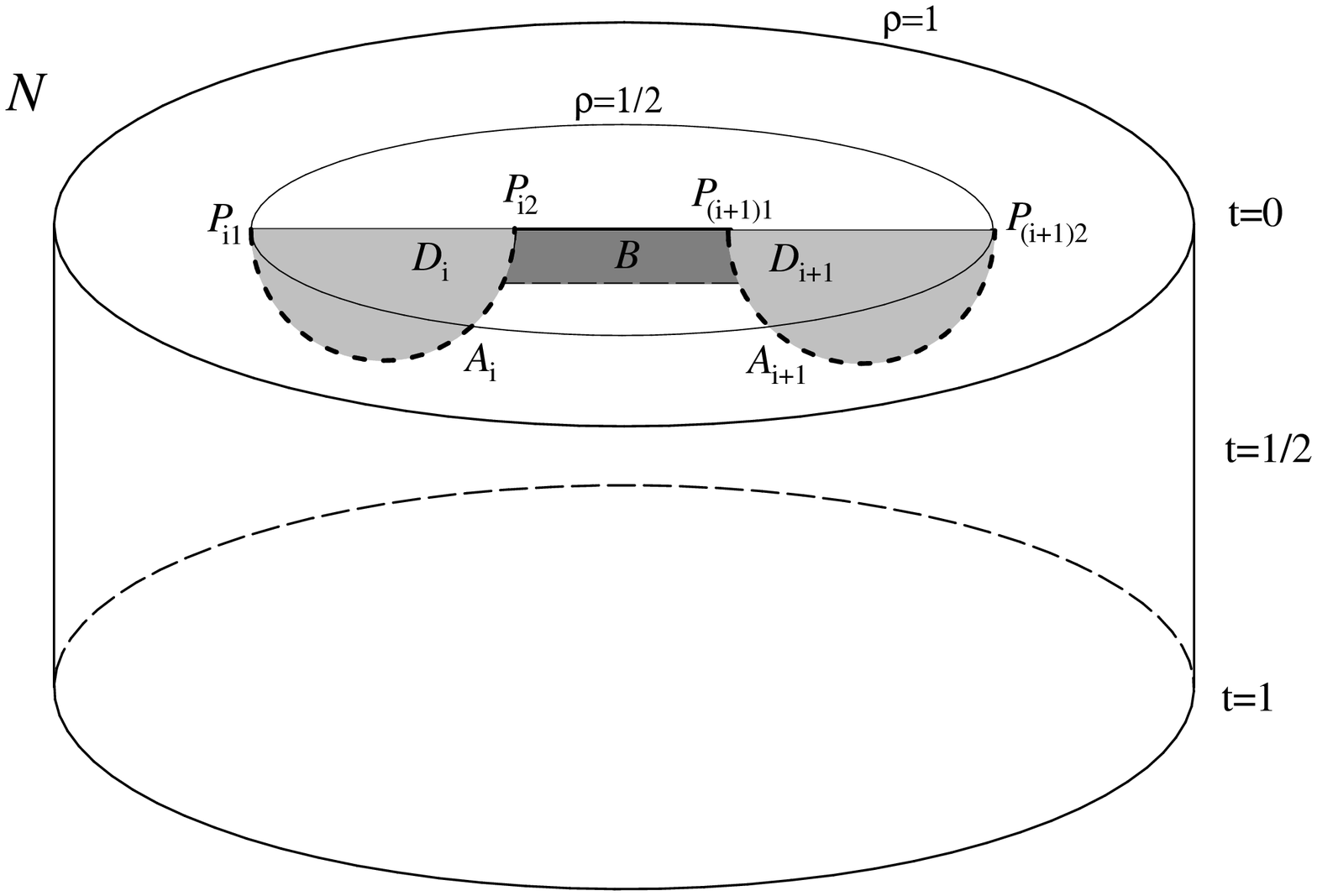}
\end{center}
\caption{} \label{finale2}
\end{figure}

 \end{description}

Moreover, let $\PMCG_{2n}(\T_g)$ be the subgroup of
$\MCG_{2n}(\T_g)$ consisting of the isotopy classes of the
homeomorphisms of $\T_g$ pointwise fixing the punctures and set
$\E_{2n}^g=\PMCG_{2n}(\T_g)\cap\X_{2n}^g$. The next proposition
shows that  $\X_{2n}^g$ is generated  by $\iota_1,\lambda_i$,
$i=1,\ldots,n-1$, and  a set of generators for $\E_{2n}^g$.

\begin{proposition}
\label{relation} Let $\Sigma_{2n}$ be the symmetric group on $2n$
letters and denote with $\Sigma'_{2n}$
 the subgroup of $\Sigma_{2n}$  generated by the transposition
$(1\ 2)$ and the permutations $(2i-1\ 2i+1)(2i\ 2i+2)$, for
$i=1,\ldots n-1$. Then the exact sequence
$1\to\PMCG_{2n}(\T_g)\to\MCG_{2n}(\T_g)\to\Sigma_{2n}\to 1$
restricts to an exact sequence $1\to\E_{2n}^g\to\X_{2n}^g \to
\Sigma'_{2n}\to 1.$
\end{proposition}
\begin{proof} First of all, we recall that the homomorphism
$p:\MCG_{2n}(\T_g)\to\Sigma_{2n}$ is obtained by considering the
permutation induced on the punctures by the elements of
$\MCG_{2n}(\T_g)$, where the puncture $P_{ij}$ corresponds to the
letter $2i+j-2$, for $i=1,\ldots, n$ and $j=1,2$. Let
$\varepsilon\in\X_{2n}^g$, since the extension of $\varepsilon$
induces a permutation of the arcs $A_1,\ldots,A_n$, it is easy to
see that $p(\varepsilon)\in\Sigma'_{2n}$. Moreover if
$p(\varepsilon)=1$ then an extension of it fixes the arcs, so
$\varepsilon\in\E_{2n}^g$. In order to complete the proof we only
need to check the surjectivity of $\X_{2n}^g \to \Sigma'_{2n}$.
This follows by observing that $\iota_1$ maps on $(1\ 2)$ and
$\lambda_i$ maps on $(2i-1\ 2i+1)(2i\ 2i+2)$, for $i=1,\ldots
n-1$.
\end{proof}

In the next two sections we will find a finite set of generators
for $\E_{2n}^g$ and $\X_{2n}^g$.

\section{The case of genus zero}
\label{zero}

We introduce the homeomorphisms whose isotopy classes generate
$\E_{2n}^0$ and describe their extensions. For
each $i=1,\ldots,n$, we set $A'_i=D_i\cap\partial \H_g$,
where $D_i$ is a trivializing disc for $A_i$ (see Figure
\ref{deco}).

\begin{description}
\item[Spin of a puncture] A simple closed oriented curve
$c\subset\partial\mathbf{B}^3$ containing a puncture $P_{j2}$ is
called \textit{admissible} if $c\cap A'_i=\emptyset$,
$i=1,\ldots,n$, $i\ne j$ and $c\cap A'_j=P_{j2}$. A spin
$s_{P_{j2},c}$ of the puncture $P_{j2}$ along an admissible curve
$c$ can be extended to $\mathbf{B}^3$ as follows. We can suppose
that a tubular neighborhood $N$ of $c$  in $\partial \mathbf{B}^3$
does not intersect the arcs $A'_i$, for $i\ne j$, and so there
exists an embedded ball $\bar N= [-1,1]\times \mathbf{D}^2$ in
$\mathbf{B}^3$ such that $[-1,1]\times\partial \mathbf{D}^2$ is
$N$, $\bar N\cap A_i=\emptyset$, for $i\ne j$ and $\bar N\cap A_j$
is an arc $l$ transversal to the trivial fibration of $\bar N$ in
discs. We can extend $s_{P_{j2},c}$ to $\mathbf{B}^3$ by the
identity outside $\bar N$ and by making the mapping cone from the
center of each disc $\{y\}\times \mathbf{D}^2$ inside $\bar N$.
Obviously such an extension keeps $A_i$ fixed for $i\ne j$. Since
two trivial properly embedded arcs in a ball with the same
endpoints are ambient isotopic by an isotopy fixing the boundary,
then $s_{P_{j2},c}$ is isotopic to a homeomorphism that keeps $l$
and so $A_j$ fixed. Therefore, it admits an extension that keeps
all the arcs fixed and, as a consequence, $ s_{P_{j2},c}\in
\E_{2n}^0$.

We denote with $s_{ji}$ the spin of the point $P_{j2}$ about the
curve depicted in Figure \ref{spin}, for $i,j=1,\ldots,n$.

\begin{figure}[ht]
\begin{center}
\includegraphics*[totalheight=5cm]{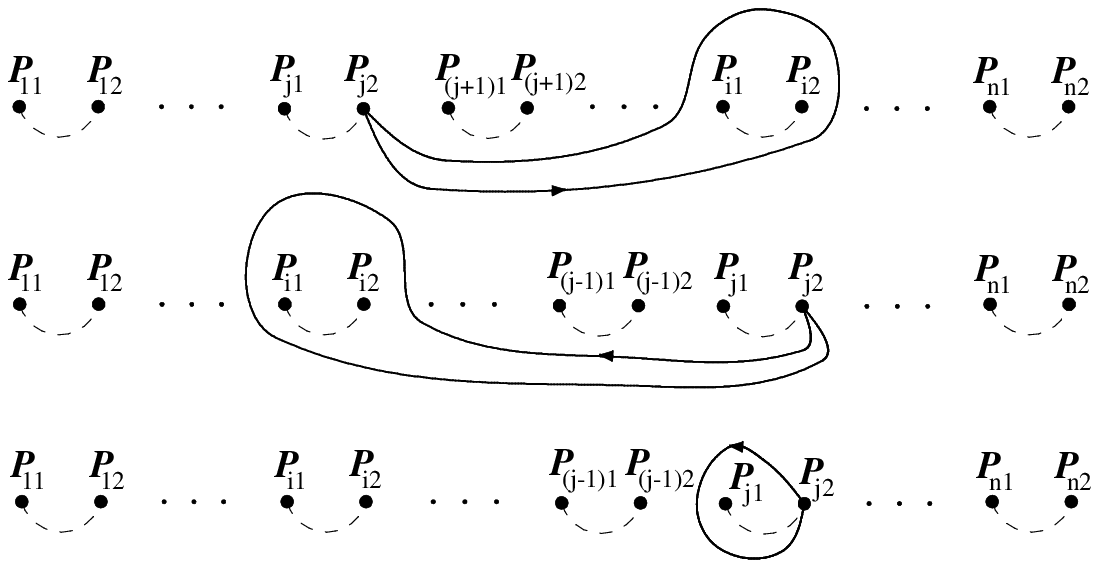}
\end{center}
\caption{The spin $s_{ji}$.} \label{spin}
\end{figure}

\begin{remark} \label{remark} Note that if $c=c_1\cdots c_k$ in
$\pi_1(\S^2-\{P_{j1}, P_{i1},P_{i2}\,\vert\,i=1,\ldots,n,\ i\ne
j\},P_{j2})$ then $s_{P_{j2},c}=s_{P_{j2},c_k}\cdots
s_{P_{j2},c_1}$ in $\PMCG_{2n}(\S^2)$. Moreover if  $c$ is
admissible, then we can take the $c_i$'s in
$\pi_1(\S^2-(\{P_{j1}\}\cup_{i\ne j} A'_i ),P_{j2})$ and so
$s_{P_{j2},c}$ is a product of the $s_{ji}$'s.
\end{remark}

\item[Slide of an arc] \label{slide} Consider a simple closed
oriented curve $c\subset\T_g-\P_n$ such that $c\cap A'_j$ is a
single point $P\ne P_{j1},P_{j2}$.  The slide $S_{j,c}$ of the arc
$A_j$ along the curve $c$ is defined as the spin of the point $P$
along the curve $c$, where we require that the parametrization of
the tubular neighborhood $N\cong[-1,1]\times\S^1$ of $c$ is chosen
such that $A'_j$ is parametrized by $[-1/2,1/2]\times \{0\}$, and
so is kept fixed by $S_{j,c}$. Note that $S_{j,c}$ is isotopic to
the product of spins $s_{P_{j1},c_1}s_{P_{j2},c_2}$, where $c_1$
and $c_2$ are the curves parametrized by $\{-1/2\}\times\S^1$ and
$\{1/2\}\times\S^1$. In order to extend $S_{j,c}$ to the handlebody, let us consider an embedded solid torus $\widetilde
N=[-1,1]\times \S^1\times [0,1]$ in $\H_g$ such that
$N=[-1,1]\times \S^1\times \{0\}$, $\widetilde N\cap
A_i=\emptyset$, for $i=1,\ldots,n$, $i\ne j$, and $\widetilde
N\cap A_j\subset [-1,1]\times \S^1\times [0,1/2]$. Then we can
extend $S_{j,c}$ to $\H_g$ by the identity outside
$\widetilde N$, and in $\widetilde N$ by the map
$$(y,\theta,t)\mapsto\left\{\begin{array}{ll}(S_{j,c}(y,\theta),t) & \textup{ if } 0\le t\le
1/2\\(y,\theta+\theta',t)& \textup{ if }1/2\le t\le
1,\end{array}\right.$$

with $$\theta'=\left\{\begin{array}{ll} -4\pi(2-2t)(y+1)+2\pi&
\textup{ if }\ -1\le y\le -1/2\\2\pi (2t-1)& \textup{ if } -1/2\le
y\le 1/2\\4\pi(2-2t)(y-1)+2\pi& \textup{ if } \ 1/2\le y\le 1
\end{array}\right..$$

It is clear that $ S_{j,c}$  fixes all the arcs and so belongs to
$\E_{2n}^g$.

We denote with $S_{ji}$ and $S'_{ji}$ the slides of the arc $A_j$
along the curves depicted in Figure \ref{track} (a) and (b),
respectively, for $i,j=1,\ldots,n$, $i\ne j$. Moreover, $L_{ki}$ and $M_{ki}$ are the slides of the arc $A_k$
along the curves  $e'_{ki}$, $g'_{ki}$ depicted in Figure \ref{uno},
for $i=1,\ldots,g$ and $k=1,\ldots,n$. 

\begin{figure}[ht]
\begin{center}
\includegraphics*[totalheight=8cm]{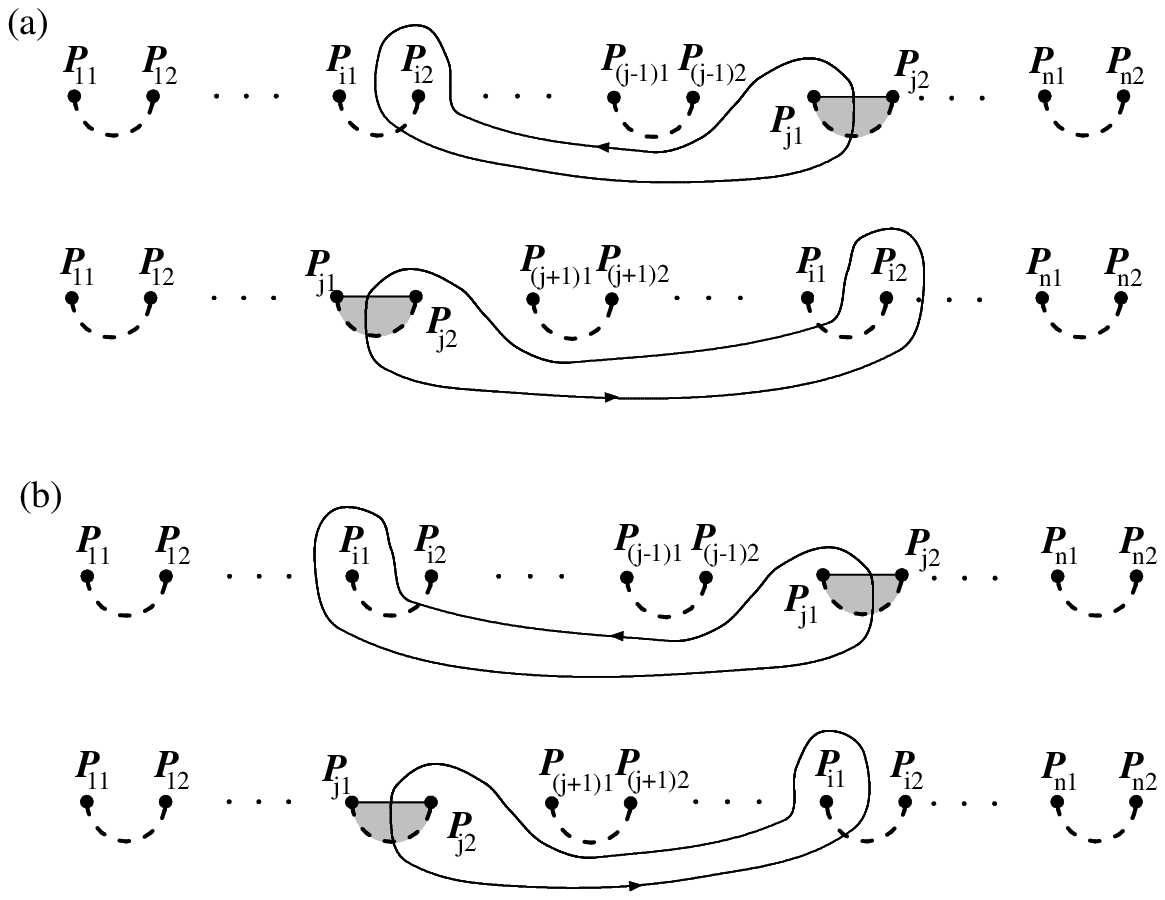}
\end{center}
\caption{Slides  $S_{ji}$ and $S'_{ji}$.} \label{track}
\end{figure}

\begin{figure}[ht]
\begin{center}
\includegraphics*[totalheight=6cm]{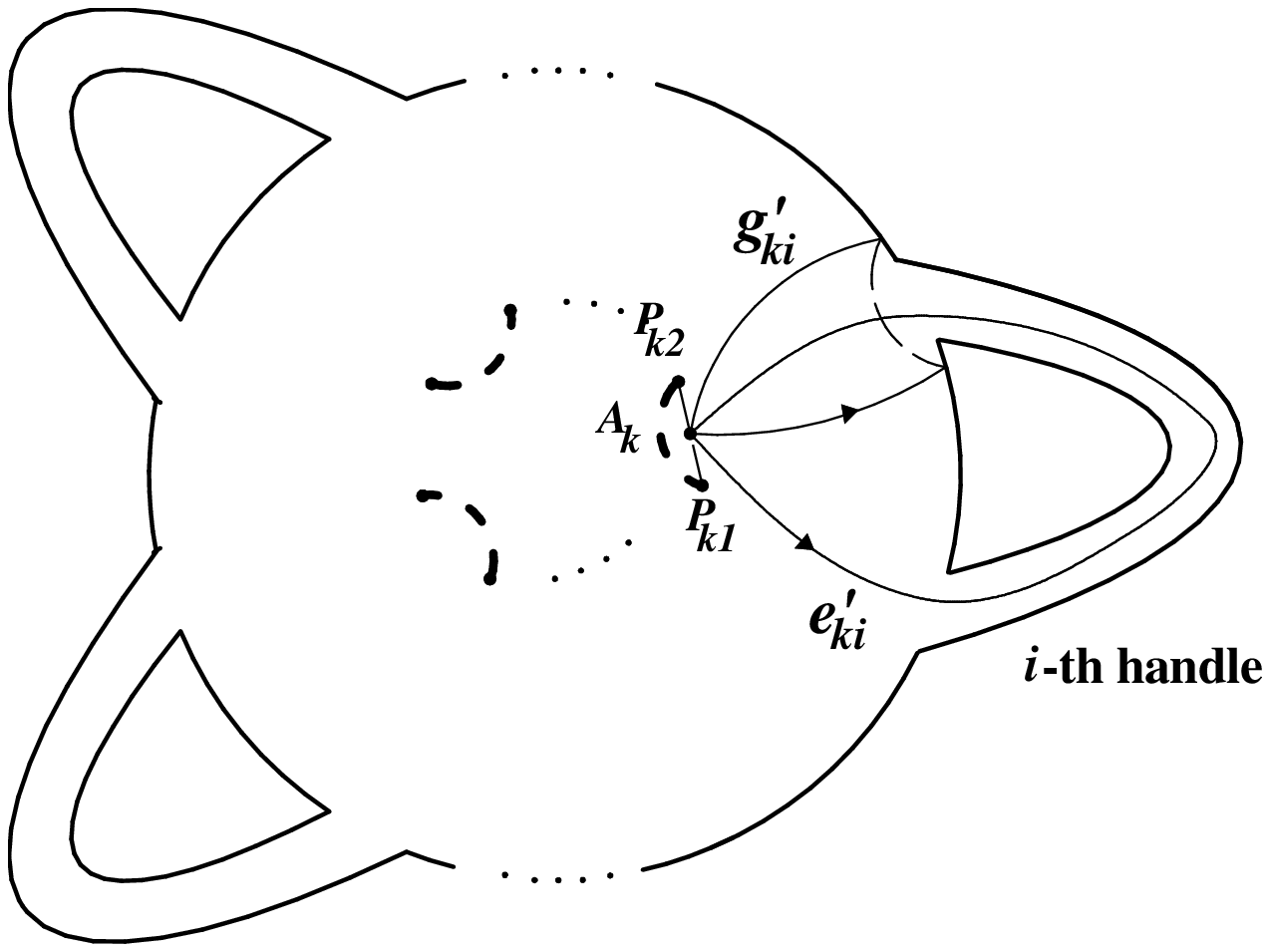}
\end{center}
\caption{Curves $e'_{ki}$ and  $g'_{ki}$.} \label{uno}
\end{figure}

\end{description}

\begin{figure}[ht]
\begin{center}
\includegraphics*[totalheight=7cm]{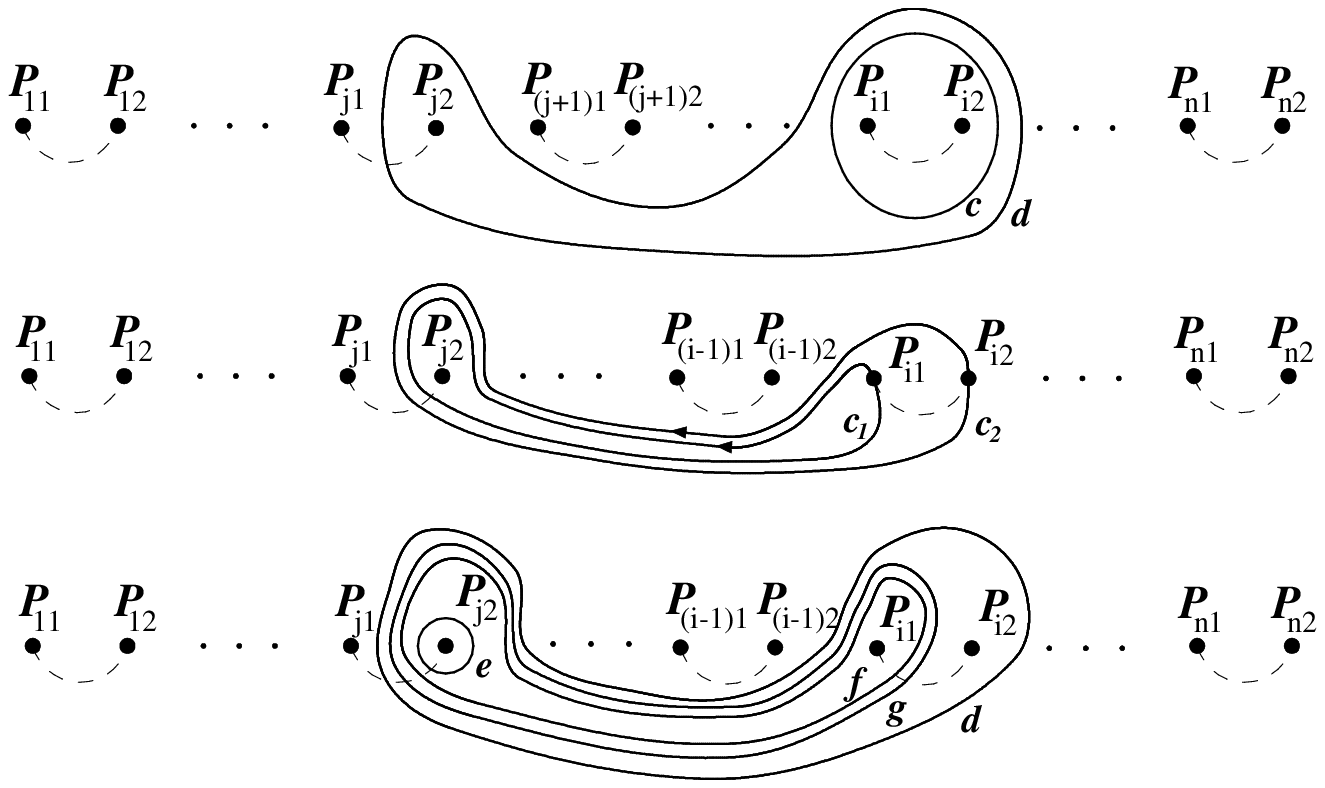}
\end{center}
\caption{} \label{ridu}
\end{figure}

\begin{lemma} \label{riduzione} If $g=0$ then the following relation holds $s_{ji}=S_{ij}s_{ii}^{-1}$
for each $i,j=1,\ldots,n$ and $i\ne j$.
\end{lemma}
\begin{proof} Let $j<i$. Referring to Figure \ref{ridu} we have that
$s_{ji}=t_{c}^{-1}t_{d}=s_{ii}^{-1}t_d$. Moreover,
$S_{ij}=s_{P_{i1},c_1}s_{P_{i2},c_2}=t^{-1}_{e}t_{f}t_{g}^{-1}t_{d}=t_d$,
where the last equality holds since $g$ and $f$ are isotopic
curves and $t_e$ is trivial. So $s_{ji}=s_{ii}^{-1}S_{ij}$. The
proof in the case $j>i$ is completely analogous.
\end{proof}

\begin{proposition}[\cite{CM}]\label{braid} Consider the homomorphism
$j_{g,n,m}:\PMCG_n(\T_g)\rightarrow\PMCG_m(\T_g)$ induced by the
inclusion, with $1\leq m<n$ for $g\geq 1$ and $3\leq m<n$ for
$g=0$. Then $\ker j_{g,n,m}\cong \pi_1(F_{m,n-m}(\T_g))$, where
$F_{m,n-m}(\T_g)$ denotes the configuration space of $n-m$ points
in $\T_g$ with $m$ punctures.
\end{proposition}

The next result describes a set of generators for the group
$\E_{2n}^0$.
\begin{theorem} \label{g=0} The subgroup $\E_{2n}^0$ of $\PMCG_{2n}(\S^2)$ is generated by
spins $s_{jj}$  of $P_{j2}$ along the simple closed oriented curve
depicted in Figure \ref{spin} and slides $S_{ji}$,  $S'_{ji}$ of
the arc $A_j$ along the curves depicted in Figure \ref{track} (a)
and (b), respectively, for $i,j=1,\ldots,n$, $i\ne j$.
\end{theorem}
\begin{proof}
Let $\mathcal{G}_{2n}^0$ be the subgroup of $\E_{2n}^0$ generated
by the following elements
\begin{itemize}
\item[a)] spins $s_{P_{j2},c}$  of $P_{j2}$ along  admissible
curves on $\S^2$, $j=1,\ldots,n$;

\item[b)] slides $S_{j,c}$ of the arc $A_j$ along  simple closed
oriented curves, in $\S^2-\P_{n}$ intersecting $A'_j$ in a single
point different from $P_{j1}$ and $P_{j2}$, $j=1,\ldots,n$.
\end{itemize}
We will  prove that $\mathcal{G}^0_{2n}=\E_{2n}^0$  by induction
on $n$, and that $\mathcal{G}^0_{2n}$ is generated by $s_{ji}$
with $i,j=1,\ldots,n$, $i\leq j$ and $S_{ji},S'_{ji}$, with
$i,j=1,\ldots,n$, $j \ne i$.

For $n=1$ the group $\E_{2}^0$ is trivial since $\PMCG_2(\S^2)$ is
trivial, so there is nothing to prove. In order to prove the
inductive step we need the following lemma.
\begin{lemma} \label{lemma}  Given an element  $h\in \E_{2n}^0$ there
exists an element $h'\in\mathcal{G}^0_{2n}$  such that $h'
h(A'_n)=A'_n$.
\end{lemma}
\begin{proof}
Let $\widetilde A_0=h(A'_n)$ and $\widetilde{D}_0=h(D_n)$. Suppose
that $\widetilde{A}_0\cap A'_i=\emptyset$ for $i=1,\ldots, n-1$.
Then the closed  curve $c=\widetilde{A}_0\cup A'_n$ based on
$P_{n2}$ is homotopic rel $P_{n2}$ to the product of simple closed
curves $c_1\cdots c_k$ in $\S^2-\{A'_i\ \mid\ i=1,\ldots n-1\}$
and each $c_j$ is a curve of type a). Moreover
$s_{P_{n2},c_k}\cdots s_{P_{n2},c_1}(\widetilde{A}_0)=A'_n$, so in
this case $h'=s_{P_{n2},c_k}\cdots s_{P_{n2},c_1}$.

Let $j\in \{1,\ldots, n-1\}$ such that $\widetilde{A}_0\cap
A'_j\ne \emptyset$. Up to a small  deformation we can suppose that
$D_j$ and $\widetilde D_0$ intersect transversally and so the
connected components of the intersection are arcs or circles. Let
us consider the circular components. By an innermost argument, it
is possible to choose one of them, let us say $C$, such that the
union of the discs bounded by $C$ on $D_j$ and $\widetilde D_0$ is
a sphere $S$ that intersects $D_j\cap \widetilde D_0$ only in $C$.
Obviously $S$ does not contain any of the arcs $A_i$ and then, by
an isotopy, the intersection $C$ can be removed. Iterating the
procedure we can remove all the circular intersections. Now, let
us consider the arcs. Since $h(A_n)=A_n$, then $h(A_n)\cap
A_j=\emptyset$. So the endpoints of the  arcs in
$D_j\cap\widetilde{D}_0$ are  points of $\widetilde{A}_0\cap
A'_j$. By an innermost argument, there exists an arc $B_0$ that
determines a disc both in $\widetilde{D}_0$ and in $D_j$, whose
union is a disc $\bar D_0$, properly embedded in the ball, that
intersects $\widetilde{D}_0\cap D_j$ only in $B_0$. Let $\Delta_0$
be the the connected component of $\mathbf{B}^3-\bar D_0$ that
does not contain $A_j$. If none of the $A_i$ for $i=1,\ldots,n-1$
is contained in $\Delta_0$ then, by isotopy, the intersection
$B_0$ can be removed.

\begin{figure}[ht]
\begin{center}
\includegraphics*[totalheight=18cm]{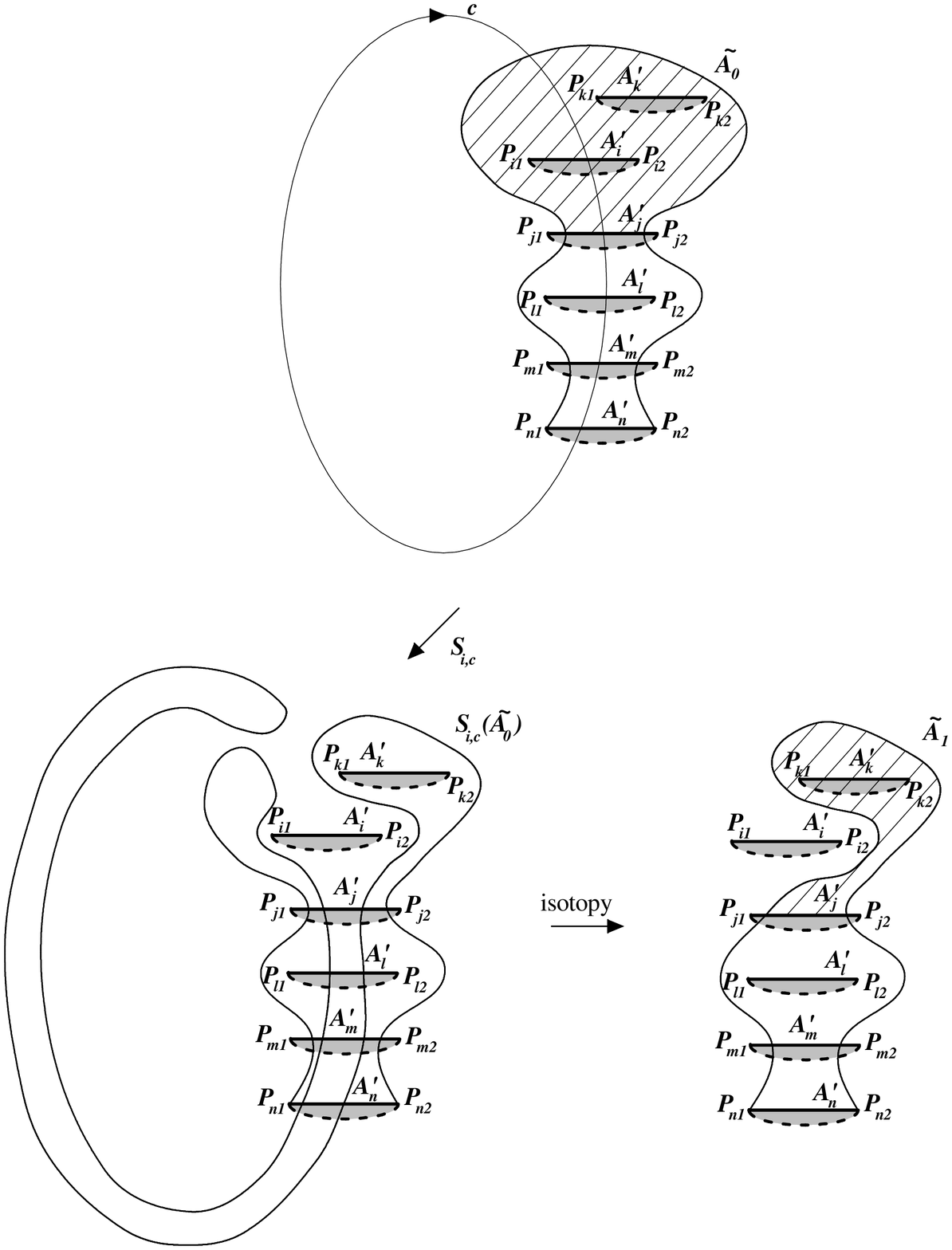}
\end{center}
\caption{} \label{proof}
\end{figure}

Otherwise, let $A_i$ be contained in   $\Delta_0$. Referring to
Figure \ref{proof}, choose a simple oriented closed curve $c$ on
$\S^2$ such that $c\cap A'_i$ is a single point and
$c\cap\widetilde{A}_0$ is a single point in $\Delta_0\cap \S^2$
(the dashed disc in the top part of the figure represents
$\Delta_0\cap\S^2$). Consider the slide $S_{i,c}$ of the arc $A_i$
along the curve $c$ and  let $h_1=S_{i,c}h$. The figure
illustrates the image of $\widetilde{A}_0$ after the application
of the homeomorphism $S_{i,c}$ and, subsequently, of a
homeomorphism isotopic to the identity. If we set
$\widetilde{A}_1=S_{i,c}(\widetilde{A}_0)=h_1(A'_n)$ and
$\widetilde{D}_1=S_{i,c}(\widetilde{D}_0)=h_1(D_n)$, we have that
$\textup{card}(\widetilde{D}_1\cap
D_j)=\textup{card}(\widetilde{D}_0\cap D_j)$, and
$\widetilde{A}_1\cap A'_k=\widetilde{A}_0\cap A'_k$, for each
$k=1,\ldots,n-1$.  Moreover, by the same argument used before, it
is possible to find an arc $B_1$ that determines a disc both in
$\widetilde{D}_1$ and in $D_j$, whose union is a disc $\bar D_1$,
properly embedded in the ball, that intersects
$\widetilde{D}_1\cap D_j$ only in $B_1$. If we denote with
$\Delta_1$ the connected component of $\mathbf{B}^3-\bar D_1$ that
does not contain $A_j$ (the dashed disc in the bottom part of the
figure represents $\Delta_1\cap\S^2$) then $\mathcal A\cap
\Delta_1=\mathcal A\cap \Delta_0-A_i$. So, after $k$-steps and
composing $h$ with the opportune slides, we obtain the case of
$\Delta_k\cap\mathcal A=\emptyset$ and so as before we can remove
the  intersection arc $B_k$. Since the intersections of
$\widetilde{D}_i\cap D_j$   are finitely many, we can remove them
all and so return to the case in which $\widetilde{A}_k\cap A'_j=
\emptyset$.

Since  during the process we do not increase the number of
intersections of $\widetilde{A}_k$ with the other arcs $A'_h$ by
applying it to the other arcs, we return to the case
$\widetilde{A}_i\cap A'_k=\emptyset$ for each $k=1,\ldots n-1$
already considered.
\end{proof}

\textit{Continuation of the proof of Proposition \ref{g=0}} By the
previous lemma it is enough to consider the subgroup
$\widetilde{\X}_{2n}^0$ of $\E_{2n}^0$, consisting of the elements
fixing $A'_n$. Let $h\in \widetilde{\X}_{2n}^0$, by hypothesis
$h(A_n)=A_n$,  so we can assume that $h$ fixes the whole disc
$D_n$, as well as a tubular neighborhood $N$ of it. By contracting
$N$ to the point $ P_{n2}$ we obtain a surjective map
$i_1:\widetilde{\X}_{2n}^0\to E^0_{2n-1}$, where $ E^0_{2n-1}$ is
the group of elements of $\PMCG_{2n-1}(\S^2)$ that extend to the
ball $\mathbf{B}^3$ fixing $A_1,\ldots,A_{n-1}$. Moreover, the
surjective homomorphism
$j_{0,2n-1,2n-2}:\PMCG_{2n-1}(\S^2)\rightarrow\PMCG_{2n-2}(\S^2)$
of Proposition \ref{braid} restricts to a surjective homomorphism
$i_2:E^0_{2n-1}\to\E_{2n-2}^0$. When $n=2$,  the kernel of $i_2$
is trivial since $\PMCG_{3}(\S^2)=\PMCG_{2}(\S^2)=1$. Otherwise,
when $n>2$, by \cite[pp. 158-160]{Bi} and \cite{CM}, $\ker i_2$ is
generated by $A=\{s_{ P_{n2},\gamma_i}\ \vert\ i=1,\ldots,2n-3\}$,
where the loops $\gamma_i$ are generators of $\pi_1(\S^2-\P_{n-1},
P_{n2})$. By the induction hypothesis, $\E_{2n-2}^0$ is generated
by $B=\{s_{P_{j2},c}\ \vert\ c\cap A'_i=\emptyset,\ c\cap
A'_j=P_{j2}, \ i,j=1,\ldots,n-1,\ i\ne j \}$ and $C=\{S_{j,c}\
\vert\ c\subset \S^2-\P_{n-1},\ \sharp(c\cap A'_j)=1,\ c\cap
A'_j\ne P_{j1},P_{j2} j=1,\ldots,n-1\}$. So a complete set of
generators for $\widetilde{\X}_{2n}^0$ is given by the generators
of $\ker i_1$, the lifts of elements in $A$ via $i_1$, and the
lifts of elements in $B$ and $C$ via $i_2i_1$. The kernel of $i_1$
is generated by $s_{nn}$ and so it is contained in
$\mathcal{G}^0_{2n}$. A spin $s_{ P_{n2},\gamma_i}\in A$ lifts to
a slide of the arc $A_n$ along a suitable  curve,
 and so belongs to $\mathcal G^0_{2n}$. An element $s_{P_{j2},c}$
in $B$ lifts to a spin still based on $P_{j2}$ along the lifting
of $c$. Moreover, since we can suppose that $c$ avoids $P_{n2}$,
its lifting does not intersect $A'_n$, and so belongs to
$\mathcal{G}^0_{2n}$. In the same way the lifting of each element
in $C$ belongs to $\mathcal{G}^0_{2n}$.

Now we find a finite set of generators for $\mathcal{G}^0_{2n}$.
From Remark \ref{remark}  and  Lemma \label{riduzione}, it follows
that elements of type a) are generated by $s_{jj}$ for and
$S_{ji}$ for $i,j=1,\ldots,n$  $i\ne j$. Note that a slide $S$ of
the arc $A_j$   fixes $A'_j$. Consider the analogue of the map
$i_2$ obtained by contracting a tubular neighborhood of $A'_j$ to
$P_{j2}$. The kernel of the map is $s_{jj}$,  while the image of
$S$ is a spin $s$ based on $P_{j2}$ along a curve in
$\S^2-\{P_{j1},P_{i1},P_{i2}\,\mid\, i=1,\ldots,n,\ i\ne j\}$.
Since $s$ decomposes into a product of spins along a set of curves
whose homotopy classes generate
$\pi_1(\S^2-\{P_{j1},P_{i1},P_{i2}\,\mid\, i=1,\ldots,n,\ i\ne
j\}, P_{j2})$
 and the lifting of each spin is one of   the $S_{ji}$'s or the $S'_{ji}$'s, then  $S$
 is a product of  the $S_{ji}$'s and $S'_{ji}$'s.
\end{proof}

\begin{corollary} \label{corollary}The subgroup $\X^0_{2n}$ of $\MCG_{2n}(\S^2)$
is generated by $\iota_1,\lambda_k,s_{11},S_{12},S'_{12}$, for
$k=1,\dots,n-1$.
\end{corollary}
\begin{proof} By Theorem \ref{braid} and Proposition
\ref{relation} the group $\X^0_{2n}$ is generated by
$\iota_1,\lambda_k,s_{jj},S_{ji},S'_{ji}$, for $j,i=1\ldots,n$,
$i\ne j$ and $k=1,\dots,n-1$. Set
$\Lambda_i=\lambda_1^{-1}\cdots\lambda_{i-1}^{-1}$, then the
statement follows from the relations
$s_{ii}=\Lambda_{i}^{-1}s_{11}\Lambda_{i}$, $S_{ji}=(\lambda_1
\Lambda_{i}\Lambda_{j})^{-1} S_{12}\lambda_1
\Lambda_{i}\Lambda_{j}$ if $j<i$ and
$S_{ji}=(\Lambda_j\Lambda_i)^{-1}S_{12}\Lambda_j\Lambda_i$ if
$j>i$.
\end{proof}

We note that the sets of generators obtained for $\X^0_{2n}$  and
$\E_{2n}^0$ is considerably smaller than the ones obtained in
\cite{H}.

\section{The general case}
\label{general}

We are now ready to analyze the general case. We introduce new 
homeomorphisms whose isotopy classes are generators of $\E_{2n}^g$ and
describe their extension to $\H_g$.

\begin{figure}[ht]
\begin{center}
\includegraphics*[totalheight=8cm]{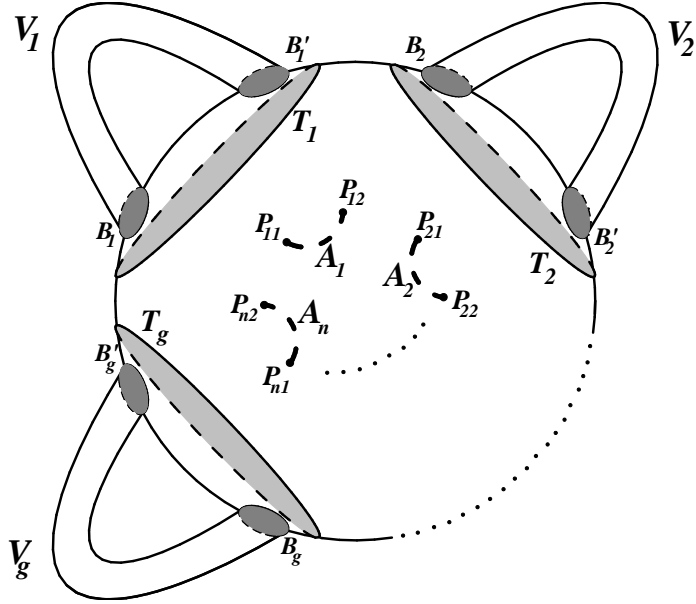}
\end{center}
\caption{A model for $\H_g$.} \label{Fig1}
\end{figure}

\begin{description}
\item[Semitwist of a handle] let $T_i$ be the 2-cell depicted in
Figure~\ref{Fig1} for $i=1,\ldots,g$. Then $T_i$ cuts away from
$\H_g$ a solid torus $K_i$, containing the $i$-th handle $\V_i$.
We denote with $\omega_i$ the homeomorphism of $\H_g$ which is  a
counterclockwise rotation of $\pi$ radians of  $K_i$ along $T_i$
on $K_i$  and the identity outside a tubular neighborhood of
$K_i$, for $i=1,\ldots, g$. As usual, we still denote with
$\omega_i$ its restriction to $\T_g$. Note that $\omega_i^2$ is
isotopic to the Dehn twist along $\partial T_i$.

\item[Twist of a meridian disk]  Referring to Figure \ref{Fig1},
let $\b_i=\partial B_i$, for $i=1,\ldots, g$. We denote with
$\tau_i$ the right-handed Dehn twist along $\b_i$. Note that
$\tau_i$ admits an extension to $\H_g$, whose effect is to give a
complete twist to the $i$-th meridian disk $B_i$.

\item[Exchanging  two handles] let $C_i$ be the properly embedded
2-cell, depicted in Figure \ref{Fig4}, for $i=1,\ldots,g-1$. Then
$C_i$ cuts away from $\H_g$ a handlebody $K'_i$ of genus two,
containing the $i$-th and the $(i+1)$-th handles. Let
$\bar{\rho}_{i}$ be the homeomorphism of $\H_g$  which exchanges
$V_i$ and $V_{i+1}$ by a counterclockwise rotation of $\pi$
radians along $C_i$ and is the identity  outside a tubular
neighborhood of $K'_i$. We set  $\rho_i=\omega_i^{-1}
\omega_{i+1}^{-1}\bar{\rho}_{i}$ for $i=1,\ldots, g-1$. Moreover,
for $i<j$ we set
$\rho_{ij}=\rho_i\rho_{i+1}\cdots\rho_{j-2}\rho_{j-1}\rho_{j-2}^{-1}\cdots\rho_{i+1}^{-1}\rho_{i}^{-1}$.
Obviously $\rho_{ij}$ exchanges $i$-th handle with the $j$-th
handle and keeps fixed the other handles.

\begin{figure}[ht]
\begin{center}
\includegraphics*[totalheight=6cm]{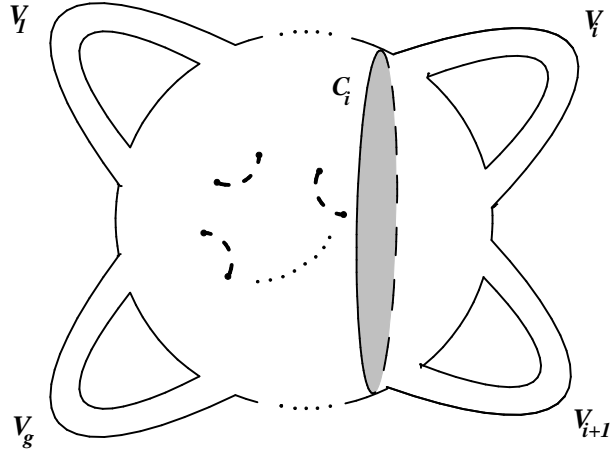}
\end{center}
\caption{Exchanging two handles.} \label{Fig4}
\end{figure}

\item [Slides of a meridian disc] Referring to Figure \ref{Fig1},
let $Z_{i}$ and $Z'_i$  be the centers of the properly embedded
meridian discs $B_i$ and $B'_i$  in $\H_g$, respectively.
Moreover, denote with $\H_g^i$ the genus $(g-1)$ handlebody
obtained  by removing the $i$-th handle $\V_i$ from $\H_g$.  A
simple closed oriented curve $e$ on $\partial\H_g^i-\P_n$
containing $Z_i$, with $B'_i \cap e=\emptyset$, will be called an
$i$\textit{-loop}. If we require that in the parametrization of a
tubular neighborhood $N\cong [-1,1]\times \S^1$ of $e$ the disk
$B_i$ is contained in $[-1/2,1/2]\times \S^1$ and $B'_i\notin N$,
then the  spin of $Z_i$ along  $e$ keeps $B_i$ and $B'_i$ fixed,
and so can be extended to $V_i$ by the identity and
 to $\H^i_g$ in the same way as the extension of a slide of an
arc (see page \pageref{slide}).  We call this homeomorphism of
$\H_g$,  as well as its restriction to $\T_g$, a \textit{slide of}
$B_i$ along the $i$-th loop $e$ and denote it with $\sigma_{i,e}$.
In a completely analogous way we can define an $i'$-loop $e'$ and
a slide of $B'_i$ along $e'$ and denote it with $\sigma'_{i',e'}$.
Note that $\sigma'_{i',e'}=\omega_i^{-1}\sigma_{i,e}\omega_i$,
where $e=\omega_i(e')$.   We set $\theta_{ij}=\sigma_{i,e_{ij}}$,
$\eta_{ij}=\sigma_{i,g_{ij}}$, $\xi_{ik}=\sigma_{i,f_{ik}}$,
$\zeta_{ik}=\sigma_{i,l_{ik}}$, where
$e_{ij},g_{ij},f_{ik},l_{ik}$ are the oriented curves depicted in
Figure \ref{Fig2}, with $i,j=1,\ldots, g$, $i\neq j$, and
$k=1,\ldots,n$.

\begin{figure}[h]
\begin{center}
\includegraphics*[totalheight=7cm]{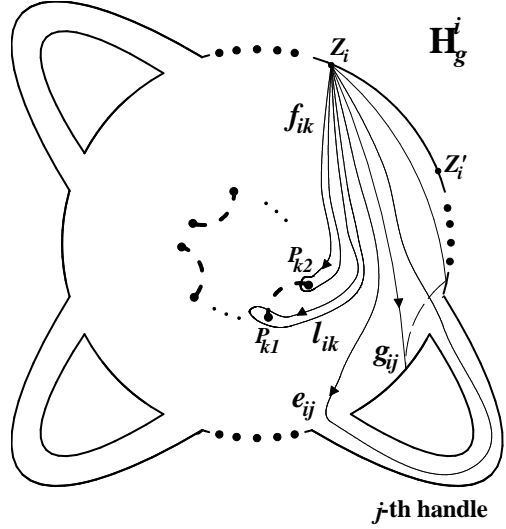}
\end{center}
\caption{The loops $e_{ij},g_{ij},f_{ik},l_{ik}$, on $\H^i_g$.}
\label{Fig2}
\end{figure}

\end{description}

\begin{remark}\label{R2} By  \cite[Lemma 3.6]{Su} if  the $i$-loop
$e$ is homotopic to the product of $i$-loops $e_1\cdots e_n$
 on $\partial \H_g^i-(\P_n\cup \{Z'_i\})$ rel
$Z_i$, then $\sigma_{i,e}$ is isotopic to
$\sigma_{i,e_n}\cdots\sigma_{i,e_1}$  modulo $\tau_i$. Since the
loops $e_{ij},g_{ij},f_{ik},l_{ik}$, for $k=1,\ldots,n$ and
$j=1,\ldots g$, with $i\ne j$, are a free set of generators for
$\pi_1(\partial \H_g^i-(\P_n\cup \{Z'_i\}), Z_i)$, then
$\theta_{ij},\eta_{ij},\xi_{ik},\zeta_{ik},\tau_i$, with
$k=1,\ldots,n$, $j=1,\ldots g$ and $i\ne j$, generate all the
slides of the $i$-th meridian disc $B_i$.
\end{remark}

Now we are ready to describe a finite set of generators for
$\E^g_{2n}$.

\begin{theorem} \label{uno}  The subgroup $\E^g_{2n}$ of $\PMCG_{2n}(\T_g)$
is generated by
$\tau_1,\omega_1,\rho_{i},\theta_{12},\eta_{12},\xi_{1k},\zeta_{1k},
s_{kk},S_{kl},S'_{kl}$ and $L_{k1}$, with $i=1,\ldots,g-1$ and $k,l=1,\ldots n$,
$k\ne l$.
\end{theorem}
\begin{proof}
Let $\mathcal{G}_{2n}^g$ be the subgroup of $\PMCG_{2n}(\T_g)$
generated by
$\tau_i,\omega_i,\rho_{i},\rho_{im},\theta_{ij},\eta_{ij},\xi_{ik},\zeta_{ik},
\sigma_{i,e},\sigma'_{i,e'} s_{kh},S_{k,c}$, where $e$ is
an $i$-loop, $e'$ is an $i'$-loop and $c$ is a curve on $\T_g$, with $i,j=1,\ldots,g$, $i\ne j$,
$i<m$ and $k,h=1,\ldots n$, $h\leq k$. 

We will  prove that $\mathcal{G}_{2n}^g=\E^g_{2n}$,  by induction
on $g$, and that $\mathcal{G}_{2n}^g$ is generated by
$\tau_1,\omega_1,\rho_{i},\theta_{12},\eta_{12},\xi_{1k},\zeta_{1k},
s_{kk},S_{kl},S'_{kl}$ and $L_{k1}$, with $i=1,\ldots,g-1$ and $k,l=1,\ldots n$,
$k\ne l$.

The case  $g=0$ is proved in Theorem \ref{g=0}.

Now we prove the inductive step. Let $g>0$ and denote with
$\widetilde{\X}_{2n}^g$  the subgroup of $\E_{2n}^g$ consisting of
the isotopy classes of the homeomorphisms that are the identity on
the boundary of the $g$-th handle $\V_g$. By the same arguments as
for  the proof of \cite[Lemma 4.4]{Su}, we have that,  for each
$h\in\E^g_{2n}$, there exists an element $h'\in
\mathcal{G}^g_{2n}$ such that $h'h$ is the identity on a meridian
disk of the $g$-th handle and so, up to isotopy, on all the
handle. Therefore $h'h\in\widetilde{\X}_{2n}^g$ and so it is
enough to show  that
$\widetilde{\X}_{2n}^g\subseteq\mathcal{G}^g_{2n}$.

Let $f:\T_g\to\T_g$ be a homeomorphism fixing $\partial \V_g$
pointwise and whose isotopy class belongs to
$\widetilde{\X}_{2n}^g$. By cutting out the $g$-th handle, and
capping  the resulting holes with the two disks $B_g$ and $B'_g$,
we can identify $f$ with a homeomorphism $f'$ of $\T_{g-1}$, such
that $f'_{|_{B_g\cup B'_g}}=\textup{Id}$ and $f'=f$ on
$\T_{g-1}-(B_g\cup B'_g)$. Moreover, by shrinking $B_g$ and $B'_g$
to their centers $Z_g$ and $Z'_g$, the map $f'$ becomes a map
$\widetilde{f}$ of $\T_{g-1}$ fixing $Z_g$ and $Z'_g$. In order to
simplify the notation, we set $P_{2n+1}=Z_g$ and $P_{2n+2}=Z'_g$.
Obviously, $\widetilde{f}$ extends to $\H_{g-1}$ fixing $\mathcal
A$ pointwise. So we obtain a surjective map
$i_1:\widetilde{\X}_{2n}^g\to E^{g-1}_{2n+2}$, where
$E^{g-1}_{2n+2}$ is the subgroup of $\PMCG_{2n+2}(\T_{g-1})$
consisting of  the elements which extend to the handlebody fixing
$\mathcal A$ pointwise. Moreover, the surjective homomorphism
$j_{g-1,2n+2,2n}:\PMCG_{2n+2}(\T_{g-1})\rightarrow\PMCG_{2n}(\T_{g-1})$
of Proposition \ref{braid} restricts to a surjective homomorphism
$i_2:E_{2n+2}^{g-1}\rightarrow\E_{2n}^{g-1}$. So, a set of
generators of  $\widetilde{\X}_{2n}^{g}$ is given by the
generators of $\ker i_1$, the lift of the generators of $\ker i_2$
via $i_1$ and the lift of the generators of $\E_{2n}^{g-1}$ via
$i_2 i_1$. The kernel of $i_1$ is generated by $\tau_g$, and so
belongs to $\mathcal{G}^g_{2n}$. By \cite[pp. 158-160]{Bi} and
\cite{CM} $\ker i_2$ is generated by  spins of $Z_g$ and $Z'_g$
about appropriate loops not containing $\P_n$, lifting to slides
of $B_g$ and $B'_g$ on $\T_g$, which are elements of
$\mathcal{G}_{2n}^g$. Moreover, by the induction hypothesis
$\E_{2n}^{g-1}=\mathcal{G}_{2n}^{g-1}$. Since we can suppose that
the generators of $\mathcal{G}_{2n}^{g-1}$ keep $Z_g$ and $Z'_g$
fixed, they lift to elements of $\mathcal{G}_{2n}^{g}$.

Now we prove that
$\tau_1,\omega_1,\rho_{i},\theta_{12},\eta_{12},\xi_{1k},\zeta_{1k},
s_{kk},S_{kl},S'_{kl}$ and $L_{k1}$, with $i=1,\ldots,g-1$ and $k,l=1,\ldots n$,
$k\ne l$ generate $\mathcal{G}^g_{2n}$. By Remark
\ref{R2}, and since, as already observed,
$\sigma'_{i',e'}=\omega_i^{-1}\sigma_{i,e}\omega_i$, where
$e=\omega_i(e')$, then the elements
$\theta_{ij},\eta_{ij},\xi_{ik},\zeta_{ik},\tau_i,\omega_i$ for
$k=1,\ldots,n$ and $i,j=1,\ldots g$, with $i\ne j$, generate all
the slides  $\sigma_{i,e}$ and $\sigma'_{i',e}$. Moreover we have
$\theta_{ij}=
\rho_{ji}^{-1}\omega_i^{-2}\theta_{ji}\omega_i^2\rho_{ji}$ if
$i>j$ and
$\theta_{ij}=\rho_{1i}\rho_{2j}\theta_{12}\rho_{2j}^{-1}\rho_{1i}^{-1}$,
if $i<j$.  The same relations hold for the other slides of a
meridian disk. 

Analogously, if $c=c_1\cdots c_j$ in
$\pi_1((\T_g-\P_{2n})\cup\{P_{k1},P_{k2}\},P)$, with $P=A'_k\cap c$,  then $S_{k,c}=S_{k,c_j}\cdots
S_{k,c_1}$, up to multiplication by $s_{kk}$. Then any slide  $S_{k,c}$ is the product of $s_{kk}$ and the slides $L_{ki},M_{ki},S_{kl},S'_{kl}$, for $i=1,\ldots,g$ and  $l,k=1,\ldots,n$ with $l\ne k$, since the corresponding cuves generate $\pi_1((\T_g-\P_{2n})\cup\{P_{k1},P_{k2}\},P)$. Moreover we have $L_{ki}=\rho_{1i}L_{k1}\rho_{1i}^{-1}$, \hbox{$M_{ki}=\rho_{1i}M_{k1}\rho_{1i}^{-1}$} and $M_{k1}=\tau_i^{-2}\zeta_{1k}\xi_{1k}s_{kk}$.
To end the proof it is enough to observe that
$\tau_i=\rho_{1i}\tau_{1}\rho_{1i}^{-1}$,
$\omega_i=\rho_{1i}\omega_{1}\rho_{1i}^{-1}$ and that, by
definition, $\rho_{ij}$ is a product of  $\rho_i$'s.
\end{proof}

We end the paper  by describing a finite set of generators for   the subgroup $\X^g_{2n}$. 
\begin{theorem}   The subgroup $\X^g_{2n}$ of $\MCG_{2n}(\T_g)$
is generated by
$\iota_1,\lambda_k,\tau_1,\omega_1,\rho_{i},\theta_{12},\eta_{12},\xi_{11},
S_{12}$ and $L_{11}$, with $i=1,\ldots,g-1$, \hbox{$k=1,\ldots,n-1$.}
\end{theorem}
\begin{proof}
The statement follows from  Proposition \ref{relation}, Corollary \ref{corollary}, Theorem \ref{uno} and the following relations: $s_{11}=\iota_1^2$, $\zeta_{1h}=\iota^{-1}_h\xi_{1h}\iota_h$, $\xi_{1h}=\lambda_{h-1}^{-1}\cdots\lambda_{1}^{-1}\xi_{11}\lambda_{1}\cdots\lambda_{h-1}$,  $S'_{12}=\iota_1^4\lambda_1^2S_{12}^{-1}$, $L_{h1}=\lambda_{h-1}^{-1}\cdots\lambda_{1}^{-1}L_{11}\lambda_{1}\cdots\lambda_{h-1}$, for each $h=1,\ldots,n$. 
\end{proof}

\bigskip\bigskip

\vspace{15 pt} {\noindent ALESSIA CATTABRIGA, Department of
Mathematics, University of Bologna, Piazza di Porta S. Donato, 5,
40126, Bologna (Italy). E-mail: cattabri@dm.unibo.it}

\vspace{15 pt} {\noindent MICHELE MULAZZANI, Department of
Mathematics, C.I.R.A.M., University of Bologna, Piazza di Porta S.
Donato, 5, 40126, Bologna (Italy). E-mail: mulazza@dm.unibo.it}


\begin{thebibliography}{5}

\bibitem{Bi}
J. S. Birman, {\it Braids, Links, and Mapping Class Groups},
Princeton University Press, Princeton-New Jersey, (1974).

\bibitem {CM1}
A. Cattabriga and  M. Mulazzani, {\it Strongly-cyclic branched
coverings of $(1,1)$-knots and cyclic presentations of groups},
Math. Proc. Cambridge Philos. Soc. {\bf 135} (2003), 137-146.

\bibitem {CM2}
A. Cattabriga and M. Mulazzani, {\it $(1,1)$-knots via the mapping
class group of the twice punctured torus}, Adv. Geom. {\bf 4}
(2004), 263-277.

\bibitem{CM} A. Cattabriga and M. Mulazzani, {\it Extending homeomorphisms from 2-punctured surfaces to
handlebodies}, to appear in Kobe J. Math., arXiv:math.GT/0601255.

\bibitem {CK}
D. H. Choi and K. H. Ko, {\it Parametrizations of 1-bridge torus
knots}, J. Knot Theory Ramifications {\bf 12} (2003), 463-491.

\bibitem {CMV}
P. Cristofori, M. Mulazzani and A. Vesnin, {\it Strongly-cyclic
branched coverings of knots via $(g,1)$-decompositions}, to appear
in Acta Math. Hungarica, arXiv: math.GT/0402393.


\bibitem {Do}
H. Doll, {\it A generalized bridge number for links in
3-manifold}, Math. Ann. {\bf 294} (1992), 701-717.



\bibitem{E}
M. Eudave-Mu$\tilde{\textup{n}}$oz, {\it Incompressible surfaces
and $(1,1)$-knots}, J. Knot Theory Ramifications {\bf 15} (2006),
935--948.




\bibitem {GHS}
H. Goda, C. Hayashi and H.-J. Song, \textit{A criterion for
satellite 1-genus 1-bridge knots}, Proc. Amer. Math. Soc.
\textbf{132} (2004),  3449-3456.



\bibitem{G}
S. Gervais, {\it A finite presentation of the mapping class group
of a punctured surface}, Topology  {\bf 40} (2001), 703-725.


\bibitem {GMM}
H. Goda, H. Matsuda and T. Morifuji, {\it Knot Floer homology of
$(1,1)$-knots}, Geom. Dedicata {\bf 112} (2005), 197-214.

\bibitem{Ha}
C. Hayashi, \textit{1-genus 1-bridge splittings for knots}, Osaka
J. Math. \textbf{41} (2004),  371-426.

\bibitem{H}
H. M. Hayashi, \textit{Generators for two groups related to the
braid group}, Pacific J. Math. \textbf{59} (1975), 475-486.

\bibitem{K}
Y. Koda, \textit{ Strongly-cyclic branched coverings and the
Alexander polynomial of knots in rational homology spheres}, to
appear in Math. Proc. Cambridge Philos. Soc.

\bibitem {LP}
C. Labru\`{e}re and L. Paris, \textit{Presentations of the
punctured mapping class groups in terms of Artin groups},  Algeb.
Geom. Topol. \textbf{1} (2001), 73-114.


\bibitem{S}
T. Saito, \textit{Genus one 1-bridge knots as viewed from the
curve complex}, Osaka J. Math. \textbf{41} (2004),  427-454.

\bibitem{Sa}
M. Sakuma,  \textit{The topology, geometry and algebra of
unknotting tunnels. Knot theory and its applications}, Chaos
Solitons Fractals \textbf{9} (1998),  739--748.



\bibitem{Su}
S. Suzuki, {\it On homeomorphisms of a 3-dimensional handlebody},
Can. J. Math. {\bf 29} (1977), 111-124.

\end{thebibliography}
\end{document}